\documentclass[letterpaper,10 pt, conference]{ieeeconf}   

\IEEEoverridecommandlockouts                              
\usepackage{graphics} 
\usepackage{epsfig} 
\usepackage{amsmath} 
\usepackage{amssymb}  
\usepackage{url}
\usepackage[latin1]{inputenc}

\newcommand{\eq}{\begin{equation}}
\newcommand{\eeq}{\end{equation}}
\newcommand{\eqn}{\begin{eqnarray}}
\newcommand{\eeqn}{\end{eqnarray}}

\newcommand{\bsea}{\begin{subeqnarray}}
\newcommand{\esea}{\end{subeqnarray}}
\newcommand{\nn}{\nonumber}

\newcommand{\de}{\mathrm{d}}
\newcommand{\tr}{\mathop{\rm tr}}  

\newcommand{\Ac}{ \mathcal{A}}

\newcommand{\Dc}{ \mathcal{D}}

\newcommand{\Qc}{ \mathcal{Q}}

\newcommand{\Cs}{ \mathbb{C}}

\newcommand{\Es}{ \mathbb{E}}

\newcommand{\Rs}{ \mathbb{R}}

\newcommand{\Zs}{ \mathbb{Z}}

\def\qed{\hfill \vrule height 7pt width 7pt depth 0pt \smallskip}

\newcounter{pippo}

\newtheorem{remark}{Remark}[section]
\newtheorem{teor}{Theorem}[section]
\newtheorem{corr}{Corollary}[section]
\newtheorem{propo}{Proposition}[section]
\newtheorem{lemm}{Lemma}[section]
\newtheorem{exam}{Example}
\newtheorem{obss}{Observation}
\newtheorem{probl}[pippo]{Problem}
\newtheorem{defn}{Definition}[section]
\newcommand{\teo}{\begin{teor}}
\newcommand{\eteo}{\end{teor}}
\newcommand{\cor}{\begin{corr}}
\newcommand{\ecor}{\end{corr}}
\newcommand{\prop}{\begin{propo}}
\newcommand{\eprop}{\end{propo}}
\newcommand{\lem}{\begin{lemm}}
\newcommand{\elem}{\end{lemm}}
\newcommand{\ex}{\begin{exam}}
\newcommand{\eex}{\end{exam}}
\newcommand{\pb}{\begin{probl}}
\newcommand{\epb}{\end{probl}}
\newcommand{\df}{\begin{defn}}
\newcommand{\edf}{\end{defn}}
\newcommand{\aprop}{\begin{apropo}}
\newcommand{\eaprop}{\end{apropo}}
\newcommand{\alem}{\begin{alemm}}
\newcommand{\ealem}{\end{alemm}}
\newcommand{\rem}{\begin{remark}}
\newcommand{\erem}{\end{remark}}
\newcommand{\oss}{\begin{obss}}
\newcommand{\eoss}{\end{obss}}

\title{\LARGE \bf Factor Analysis of Moving Average Processes}

\author{Mattia~Zorzi, Rodolphe~Sepulchre\thanks{This paper presents research
results of the Belgian Network DYSCO (Dynamical Systems, Control, and
Optimization), funded by the Interuniversity Attraction Poles Programme, initiated
by the Belgian State, Science Policy Office. The scientific responsibility
rests with its authors. This research is also supported by FNRS (Belgian Fund for Scientific Research).} \thanks{M. Zorzi is is with the
Dipartimento di Ingegneria dell'Informazione, Universit\`a di
Padova, via Gradenigo 6/B, 35131 Padova, Italy,
({\tt\small zorzimat@dei.unipd.it})} \thanks{R. Sepulchre is with the
Department of Engineering, University of Cambridge, Cambridge CB2 1PZ, UK,
({\tt\small r.sepulchre@eng.cam.ac.uk}), and the
Department of Electrical Engineering and Computer Science,
University of Liège, 4000 Liège, BE.} }

\begin{document}

\maketitle
\thispagestyle{empty}
\pagestyle{empty}

\begin{abstract}
The paper considers an extension of factor analysis to moving average processes. The problem is formulated as  a rank minimization of a suitable spectral density. It is shown that it can be efficiently approximated via a trace norm convex relaxation.\end{abstract}

\section{Introduction}
Factor models  are used to compress information contained in a high dimensional data vector into a small number of common factors. Those common factors represent nonobserved variables influencing  the observations. Such models have been initially developed by psychologists for statistical tests of mental abilities, \cite{spearman_1904,BURT_1909,THURSTONE_1934} and successively in econometrics and control engineering,
\cite{LEDERMANN_1937,KALMAN_SELECTION_ECONOMETRICS_1983,SCHUPPEN_1986,PICCI_1987}.

The ``standard'' factor model is a zero mean Gaussian $n$-dimensional random vector whose covariance matrix $X$ can be decomposed as the sum of a low rank covariance matrix $Y$ plus a
diagonal covariance matrix $Z$, i.e. $X=Y+Z$. Here $X$, $Y$ and $Z$ belong to the vector space of symmetric matrices of dimension $n$, say $\mathbf{Q}_n$,
and are positive semidefinite, say $X,Y,Z\succeq 0$. $X$ encodes the observed information from data, $Y$ the compressed information through $r:=\mathrm{rank}(Y)$ independent common factors, and $Z$ the information which cannot be compressed. The estimation of $Y$ from $X$, leads to a constrained rank minimization of $Y$, a nonconvex problem and computationally NP hard. It admits a tight convex relaxation
which minimizes the trace of $Y$ and was introduced in factor analysis by Jackson and Agunwamba, \cite{JACKSON_LOWER_BOUNDS_1977}, and independently Bentler and Woodward, \cite{BENTLER_1980}.
Likewise, convex relaxations of rank minimization problems have been a topic of active research in the recent years \cite{FAZEL_MIN_RANK_APPLICATIONS_2002,FAZEL_MINIMUM_RANK_2004,EXACT_MATRIXCOMPL_CANDES,GUARANTEED_M_RANK_2010}.

In the above formulation data are modeled as originating from independent, identically distributed Gaussian random vectors and factor analysis consists only in dimension reduction in the cross-sectional dimension (i.e. the number of observed variables). A generalization is to assume that data originate from a stochastic process,
 thus compressing information in the cross-sectional and in the time dimension, \cite{GEWEKE_DYNAMIC,PICCI_PINZONI_1986,Pena_BOX__1987,ON_THE_PENA_BOX_MODEL_2001}.
Different approaches have been considered to tackle the corresponding minimum rank problem. In
\cite{GEWEKE_DYNAMIC}, Geweke considers a stationary process and approximates its spectral density by a piecewise constant function.
The factor analysis is then performed piecewise by adapting the maximum-likelihood (ML) estimation method for Gaussian random vectors. Alternatively, in \cite{ENGLE_ONE_FACTOR_1981,WATSON_ALGORITHMS_1983,Pena_BOX__1987}, the authors consider a special dynamic factor model wherein the common (dynamic) factors are only combined in a static way.

A moving average (MA) Gaussian process is obtained by filtering white Gaussian noise with an all-zero filter. Although the estimation of an MA processes is relatively simple, it is not clear how to
extract the compressible information from it. The present paper considers the case of a dynamic factor model generated by an MA process.
We show that the convex relaxation for the static case then admits a natural generalization wherein the covariance matrix $X$ is replaced by the spectral density of the MA process.
In Section \ref{section_Factor_models} we recall the standard factor analysis. In Section \ref{section_ID_procedure} we introduce its dynamic MA generalization. In Section \ref{section_min_trace} we analyze
the constrained convex optimization problem which relaxes the corresponding minimum rank problem, and in Section \ref{section_algo} we propose a matricial SDP algorithm for computing a solution to the problem. Finally, in Section \ref{section_sim} we present some simulation studies.

Throughout the paper we use the following notation.
 Functions defined on the unit circle are denoted by capital Greek letters, e.g. $\Psi(e^{i\vartheta})$ with $\vartheta\in[-\pi,\pi]$.  If $\Psi$ is positive semidefinite on the unit circle we write $\Psi \succeq 0$.
$\Ac_n$ is the linear space of $\Cs^{n \times n}$-valued  analytic
functions defined on the unit circle. The (normal) rank of $\Psi\in \Ac_n$ is defined as
\eq \mathrm{rank}(\Psi)=\underset{\vartheta\in[-\pi,\pi]}{\max} \mathrm{rank}(\Psi(e^{i\vartheta})).\eeq
We define the norm
\eq \|\Psi\|:=\underset{\vartheta\in[-\pi,\pi]}{\max}\; \sigma(\Psi(e^{i\vartheta}))\eeq
where $\sigma(X)$ denotes the maximum singular value of the matrix $X$, which is equal to its maximum eigenvalue when $X\succeq 0$.

\section{Standard Factor Analysis}\label{section_Factor_models}
The standard {\em factor model} is a static linear model
\eq \label{DFA_static_model} x=Aw_y+Bw_z\eeq
where $A\in\Rs^{n\times r}$ with $r\ll n$, $B\in\Rs^{n\times n}$ diagonal. $w_y$ and $w_z$
are Gaussian random vectors with zero mean and covariance matrix equal to the identity of dimension $r$ and $n$,
respectively. Moreover, $w_y$ and $w_z$ are independent, i.e. $\Es[w_y w_z^T]=0$. Let $y:=A w_y$ and $z:= B w_z$. The $n$-dimensional random vector
$x$ is called observed vector because some statistics of it are available. To explain the reason why (\ref{DFA_static_model}) is referred to as factor model, we need some further notation. Let $x=\left[
         \begin{array}{ccc}
           x_1 & \ldots & x_n \\
         \end{array}
       \right]^T
$, $w_y=\left[
         \begin{array}{ccc}
           w_{y,1} & \ldots & w_{y,r} \\
         \end{array}
       \right]^T
$, $w_z=\left[
         \begin{array}{ccc}
           w_{z,1} & \ldots & w_{z,n} \\
         \end{array}
       \right]^T
$, $a_{jk}$ denote the entry in position $(j,k)$ of the matrix $A$, and $b_j$ denote the $j$-th entry in the main diagonal of $B$.
Therefore, \eq x_j=\sum_{k=1}^r a_{jk} w_{y,k}+b_j w_{z,j}\eeq
namely the $j$-th observed variable $x_j$ is generated by $r$ independednt common factors $w_{y,1},\ldots, w_{y,r}$ and by the specific factor $w_{z,j}$.

In view of (\ref{DFA_static_model}), $x$ is a Gaussian random vector with zero mean and covariance matrix
denoted by $X$.
Since $w_y$ and $w_z$ are independent, we get
\eq \label{STATIC_decomp}X=Y+Z\eeq where $Y=AA^T$, with rank equal to $r$, and $Z=BB^T$, diagonal, are the covariance matrices of $y$ and $z$, respectively.
Therefore, the covariance matrix of $x$ is the sum of a low rank covariance matrix, describing the common factors, and a diagonal covariance matrix, describing the specific factors.

The purpose of factors analysis consists in characterizing common factors, representing the compressed information, and specific factors from some statistics of the observed vector.
This problem can be formalized as follows: find the decomposition low rank plus diagonal (\ref{STATIC_decomp}) from an estimate of $X$.
A natural strategy for factor analysis
is to solve the following minimum rank problem
\eqn \label{MIN_RANK_PB_statico} &\underset{Y,Z\in\mathbf{Q}_n}{\min}& \mathrm{rank}(Y)\nn\\
   &\hbox{ subject to }&   Y,Z \succeq 0 \nn\\
   && Z \hbox{ diagonal}\nn\\
  &&   X=Y+Z.\eeqn
This problem, however, is computationally NP-hard. Then, the following minimum trace problem has been proposed as approximation of (\ref{MIN_RANK_PB_statico})
\eqn  &\underset{Y,Z\in\mathbf{Q}_n}{\min}& \tr(Y)\nn\\
   &\hbox{ subject to }&   Y,Z \succeq 0 \nn\\
   && Z \hbox{ diagonal}\nn\\
  &&   X=Y+Z.\eeqn
 It turns out that the relaxed problem recovers the correct decomposition
under weak assumptions. Moreover, its solution is unique, \cite{MIN_RANK_SHAPIRO_1982}.
The reason why the approximation is effective is because the {\em convex hull} of $\mathrm{rank}(Y)$ over the set $\{ Y\in\mathbf{Q}_n \hbox{ s.t. } Y\succeq 0, \| Y\|_2\leq 1 \}$ is the trace of $Y$, here $\|Y\|_2$ denotes the spectral norm of $Y$, see \cite{FAZEL_MIN_RANK_APPLICATIONS_2002}.

\section{Moving Average Factors Analysis}\label{section_ID_procedure}
In this paper we consider the following dynamic factor model
\eqn \label{DFA_model} x(t)=\sum_{k=0}^m A_k w_y(t-k)+\sum_{k=0}^m B_k w_z(t-k),\;\; t\in\Zs
\eeqn where $A_k\in\Rs^{n \times r}$ with $r\ll n$, $B_k\in\Rs^{n \times n}$ diagonal, $w_y$ and $w_z$ are $r$ and $n$-dimensional white Gaussian noise with zero
mean and variance equal to the identity, respectively, and such that \eq \label{ortog_prop}\Es[w_y(t)w_z(s)^T]=0\;\; \forall\; t,s. \eeq Note that, (\ref{DFA_model}) is
a linear combination of two white Gaussian noises, whereas a standard MA process involves just one noise term. On the other hand, only with (\ref{DFA_model}) one characterizes the compressible information, and thus the common factors, in the data.
Similarly to the static case, we define the stochastic processes
\eq \label{def_x}y(t):=\sum_{k=0}^m A_k w_y(t-k) \eeq
\eq  \label{def_e} z(t):=\sum_{k=0}^m B_k w_z(t-k). \eeq

Model (\ref{DFA_model}) is completely described by its spectral density
\eq \Psi_x(e^{i\vartheta})=\sum_{k=-m}^m  e^{-ik\vartheta} R_k  \eeq
where $R_k:=\Es[x(t)x(t+k)^T]$ is the $k$-th covariance lag of $x$. In view of (\ref{ortog_prop}), we get
\eq \label{DFA_freq_repr} \Psi_x(e^{i\vartheta})=\Psi_y(e^{i\vartheta})+\Psi_z(e^{\vartheta})\eeq
where $\Psi_y$ is the spectral density of $y$ and $\Psi_z$ the one of $z$. Moreover, from (\ref{def_x}) and (\ref{def_e}) we obtain
\eqn  \Psi_y(e^{i\vartheta})&=&  \left( \sum_{k=0}^{m} e^{-i k\vartheta} A_k\right)\left( \sum_{k=0}^{m} e^{-i k\vartheta} A_k\right)^* \nn\\
 \Psi_z(e^{i\vartheta})&=& \left( \sum_{k=0}^{m} e^{-i k\vartheta} B_k\right)\left( \sum_{k=0}^{m} e^{-i k\vartheta} B_k\right)^*. \eeqn
Therefore, $\Psi_x$, $\Psi_y$ and $\Psi_z$ belong to the following family of pseudo-polynomial matrices in $e^{i\vartheta}$:
\eq \Qc_{n,m}=\left\{ \sum_{k=-m}^m e^{-i k\vartheta} Q_k ,\; Q_k=Q_{-k}^T\in \Rs^{n \times n} \right\}.\eeq
Moreover spectral densities must be positive semidefinite on the unit circle, hence
$\Psi_x,\Psi_y,\Psi_z\succeq 0$. Since $A_k\in\Rs^{n\times r}$ and $B_k$ is diagonal, $\Psi_y$ has (normal) rank $r$ and $\Psi_z$ is diagonal.
We conclude that the observed process $x$ of the MA factor model (\ref{DFA_model}) has spectral density $\Psi_x$ which
is given by the sum of a low rank and a diagonal matrix function belonging to $\Qc_{n,m}$ and positive semidefinite on the unit circle.

Factor analysis of the model (\ref{DFA_model}) can be formulated as follows.\smallskip \smallskip
\pb Let $\mathrm{x}(1),\mathrm{x}(2),\dots, \mathrm{x}(N)$ be a finite-length sequence extracted from a realization of $x$ and assume that $m$ is given. Find the decomposition low rank plus diagonal (\ref{DFA_freq_repr})
from $\mathrm{x}(1),\mathrm{x}(2),\dots, \mathrm{x}(N)$.\epb\smallskip \smallskip
In the above problem we assumed to know $m$. If not, one can estimate $m$ from the data by using model order selection techniques, see for instance \cite{GOODWIN_1992}.

We propose the following identification procedure for finding (\ref{DFA_freq_repr}):
\begin{enumerate}
  \item Estimate $\Psi_x(e^{i\vartheta})=\sum_{k=-m}^m e^{-ik \vartheta } R_k$, such that $\Psi_x \succeq 0$, from $\mathrm{x}(1),\mathrm{x}(2),\dots, \mathrm{x}(N)$;
  \item Compute $\Psi_y$ and $\Psi_z$ by solving the following minimum rank problem
\eqn \label{MIN_RANK_PB} &\underset{\Psi_y,\Psi_z\in\Qc_{n,m}}{\min}& \mathrm{rank}(\Psi_y)\nn\\
   &\hbox{ subject to }&   \Psi_y, \Psi_z \succeq 0 \nn\\
   && \Psi_z \hbox{ diagonal}\nn\\
  &&   \Psi_x=\Psi_y+\Psi_z .\eeqn
\end{enumerate} {\em Step 1.} This is an MA parameter estimation problem. One would compute the correlogram of $x$ and then truncate it with a $m$-length rectangular window  according to the {\em Blackman-Tukey} method, \cite[page 38]{STOICA_MOSES_SPECTRAL_1997}. However, the truncated estimate may fail to be positive semidefinite over the unit circle, especially when $m\ll N$, that is it is not a spectral density. One could overcome this problem
designing a window which preserves the positivity of the windowed correlogram. The design of this window depends
on the specific application. Since this is not the main issue we address in the paper, we consider the {\em Durbin's method}, \cite{MA_ESTIMATION_BURBIN_1959}. The sketch of this procedure is as follows:
 \begin{itemize}
   \item Fit an autoregressive (AR) model of order $\tilde m=2 m$ from $\mathrm{x}(1),\mathrm{x}(2),\dots, \mathrm{x}(N)$;
   \item Approximate the AR model with an MA model, of order $m$, via the least square method;
   \item Let $x(t)=\sum_{k=0}^{ m} C_k e(t-k)$ be the estimated MA model where $e$ is white Gaussian noise with zero mean and covariance matrix equal to the identity, then \eq \Psi_x(e^{i\vartheta })=\left(\sum_{k=0}^{ m} C_k e^{-i k\vartheta}\right)\left(\sum_{k=0}^{m} C_k e^{-i k\vartheta}\right)^*.\eeq
 \end{itemize}

{\em Step 2.} The minimum rank Problem (\ref{MIN_RANK_PB}) is computationally NP-hard. This lead us to relax it in a such way to obtain a tractable convex optimization-based method.
More precisely, we would like to approximate the rank function with a convex function. Next section is devoted to this task.

\smallskip \smallskip \rem We assumed that the MA processes (\ref{def_x}) and (\ref{def_e}) have the same order, that is $m$, for simplicity. Let $\Psi_x\in\Qc_{n,m_x}$, $\Psi_y\in\Qc_{n,m_y}$ and $\Psi_z\in\Qc_{n,m_z}$.
Clearly, $m_x$ is set in Step 1 in the above procedure. If we choose $m_x> \min \{m_y,m_z\}$, the semi-definite decomposition $\Psi_x=\Psi_y+\Psi_z$, with $\Psi_z$ diagonal, may not exist or the solution for $\Psi_y$ may be trivially full rank. If we choose $m_x= \min \{m_y,m_z\}$ such decomposition does exists, but it implies that $\Psi_y$ and $\Psi_z$
belong to $\Qc_{n,m_x}$. The unique interesting case is $m_x <\min \{m_y,m_z\}$. Note that, it is not difficult to construct examples where the low rank plus diagonal semi-definite decomposition $\Psi_x=\Psi_y+\Psi_z$ is such that $m_x<\min\{m_y,m_z\}$. Without loss of generality assume that $m_x<m_y\leq m_z$. Then, $\Psi_x$ and $\Psi_y$ can be understood as elements in $\Qc_{n,m_z}$,
setting the last $m_z-m_x$ lags of $\Psi_x$ equal to zero and with the linear constraint that the last $m_z-m_y$ lags of $\Psi_y$ are equal to zero.
Accordingly, the results we will present can be easily adapted to this case. It is worth noting condition $m_x<\min\{m_y,m_z\}$ means that we permit an ``expansion'' in the time dimension, with respect to $x$, to allow a more effective compression in the cross-sectional dimension.\erem

\section{Relaxed minimum rank problem}\label{section_min_trace}
By replacing $\mathrm{rank(\Psi_y)}$ in (\ref{MIN_RANK_PB}) with a convex function we obtain a constrained convex optimization problem.
The tightest convex lower approximation of a nonconvex function is defined as follows.
 \smallskip \smallskip  \df Given $f:\Dc \rightarrow [-\infty,\infty] $, the {\em convex hull} $\mathrm{co}\; f$ is defined as the greatest convex function such that \eq \mathrm{co}\; f(x) \leq f(x),\;\; \forall x\in \Dc.\eeq\edf \smallskip\smallskip
In \cite{ZORZI_SEPULCHRE_2014}, we prove the following result.
\smallskip \smallskip
\prop \label{prop_cov_envelo2pe} Let $\Psi_y\in\Ac_n$ be an arbitrary analytic function such that $\Psi_y\succeq 0$, we define the following restricted rank function
\eq \mathrm{rank}^\prime(\Psi_y):=\left\{
  \begin{array}{ll}
    \mathrm{rank}(\Psi_y), & \|\Psi_y\|\leq 1 \\
    +\infty, & \hbox{otherwise.}
  \end{array}
\right.\eeq
Then, the convex hull of $\mathrm{rank}^\prime(\Psi_y)$ is \eq \label{def_convex_hull}\mathrm{co}\; \mathrm{rank}^\prime(\Psi_y):=\left\{
  \begin{array}{ll}
    \tr\int_{-\pi}^{\pi} \Psi_y(e^{i\vartheta}) \frac{\mathrm{d}\vartheta}{2\pi}, & \|\Psi_y\|\leq 1 \\
    +\infty, & \hbox{otherwise}.
  \end{array}
\right. \eeq\eprop \smallskip \smallskip

Consider the following convex optimization problem:
\eqn \label{minT_PB}&\underset{\Psi_y,\Psi_z\in\Qc_{n,m}}{\min}& \tr\int_{-\pi}^{\pi} \Psi_y(e^{i\vartheta}) \frac{\mathrm{d}\vartheta}{2\pi}\nn\\
   &\hbox{ subject to }&   \Psi_y,\Psi_z \succeq 0 \nn\\
   && \Psi_z \hbox{ diagonal}\nn\\
  &&   \Psi_x=\Psi_y+\Psi_z .\eeqn

\smallskip \smallskip  \prop \label{prop_esistenza}Assume that $\Psi_x\in\Qc_{n,m}$ is positive semidefinite and bounded on the unit circle. Then Problem (\ref{minT_PB}) does admit solution.\eprop \smallskip \smallskip
\proof In (\ref{minT_PB}), $\Psi_z=\Psi_x-\Psi_y$ does not appear in the objective function, thus the optimization problem is equivalent
to \eqn  \label{minT_PB_modificato}&\underset{\Psi_y\in\Qc_{n,m}}{\min}& \tr\int_{-\pi}^{\pi} \Psi_y(e^{i\vartheta}) \frac{\mathrm{d}\vartheta}{2\pi}\nn\\
   &\hbox{ subject to }&  0 \preceq \Psi_y \preceq \Psi_x \nn\\
   && \Psi_x-\Psi_y \hbox{ diagonal}. \eeqn
Hence, it is sufficient to show that (\ref{minT_PB_modificato}) admits solution for proving the statement.

The point $\Psi_y=\Psi_x$ is feasible for Problem (\ref{minT_PB_modificato}) because $\Psi_x\succeq 0$ by assumption. Thus, the feasibility set
\eqn && \hspace{-0.8cm}\mathbf{K}=\left\{\Psi_y\in \Qc_{n,m}\hbox{ s.t. }  0\preceq \Psi_y \preceq \Psi_x,\; \Psi_x-\Psi_y \hbox{ diagonal} \right\}\nn\eeqn
is nonempty. Moreover, $\mathbf{K}$ is bounded, closed and contained in the finite dimensional space $\Qc_{n,m}$. Accordingly,  $\mathbf{K}$
is a compact set. Since the objective function in (\ref{minT_PB_modificato}) is a continuous function, by {\em Weierstrass}' theorem Problem (\ref{minT_PB_modificato}) admits a minimum over $\mathbf{K}$.\qed\\

By Proposition \ref{prop_esistenza}, Problem (\ref{minT_PB}) with $\Psi_x$ estimated as explained in Section \ref{section_ID_procedure} admits a solution. Define
$c:=\|\Psi_x\|$. Then, Problem  (\ref{minT_PB}) is equivalent to
\eqn &\underset{\Psi_y\in\Qc_{n,m}}{\min}& \frac{1}{c}\tr\int_{-\pi}^{\pi} \Psi_y(e^{i\vartheta}) \frac{\mathrm{d}\vartheta}{2\pi}\nn\\
   &\hbox{ subject to }&  0 \preceq \Psi_y \preceq \Psi_x \nn\\
   && \Psi_x-\Psi_y \hbox{ diagonal}.\eeqn
Note that the feasibility set is contained in $\tilde{\mathbf{K}}:=\{\Psi_y\in\Qc_{n,m} \hbox{ s.t. } \|\Psi_y\|\leq c \}$ and $\frac{1}{c}\tr\int_{-\pi}^{\pi} \Psi_y(e^{i\vartheta}) \frac{\mathrm{d}\vartheta}{2\pi}$ is the convex hull of $\mathrm{rank} (\Psi_y)$
over $\tilde{\mathbf{K}}$. We conclude that (\ref{minT_PB}) is the convex relaxation of the minimum rank Problem (\ref{MIN_RANK_PB}).

 \section{A matricial SDP algorithm}\label{section_algo}

 The computation of a solution to Problem (\ref{minT_PB}) requires a matrix parametrization of the problem. To this end, we
 consider $Y\in\mathbf{Q}_{n(m+1)}$ partitioned as follows
  \eq Y=\left[
              \begin{array}{cccc}
                Y_{00} & Y_{01} & \ldots & Y_{0m} \\
                Y_{10} & Y_{11} &  & \vdots \\
                \vdots &  & \ddots &  \\
                  Y_{m0} & \ldots &  & Y_{mm} \\
              \end{array}
 \right]\eeq and define the {\em shift operator}
 \eq \Delta(e^{i\vartheta})=\left[
            \begin{array}{cccc}
              I_n & e^{i\vartheta} I_n &\ldots  &e^{i m\vartheta}I_n \\
            \end{array}
          \right].
 \eeq Then
 \eq \label{representation_matricial}\Delta(e^{i\vartheta}) Y \Delta(e^{i\vartheta})^*=\mathrm{D}_0(Y)+\sum_{k=1}^m \mathrm{D}_k(Y)e^{-ik\vartheta}+\mathrm{D}_k(Y)^T e^{ik \vartheta}\eeq
 where
 \eqn \mathrm{D}_0 &:& \mathbf{Q}_{n(m+1)}\rightarrow \mathbf{Q}_n \nn\\
                  && Y\mapsto \sum_{j=0}^{m} Y_{jj}\nn\\
  \mathrm{D}_k &:& \mathbf{Q}_{n(m+1)}\rightarrow \Rs^{n\times n} \nn\\
 && Y\mapsto \sum_{j=0}^{m-k} Y_{j, j+k} \eeqn
where $k=1\ldots m$.  Therefore, $\Delta(e^{i\vartheta})Y \Delta(e^{i\vartheta})^*\in\Qc_{n,m}$.
Moreover, any element in $\Qc_{n,m}$  admits the representation (\ref{representation_matricial})
because $\mathrm{D}_k$s are surjective maps and $\mathrm{D}_j(Y),\mathrm{D}_k(Y)$ with $j\neq k$ depend on different subblocks of $Y$.

This lead us to parameterize $\Psi_y,\Psi_z\in\Qc_{n,m}$ as
\eqn \Psi_y(e^{i\vartheta})&=&\Delta(e^{i\vartheta}) Y \Delta(e^{i\vartheta})^*\nn\\
\Psi_z(e^{i\vartheta})&=&\Delta(e^{i\vartheta}) Z \Delta(e^{i\vartheta})^*\eeqn
with $Y,Z\in \mathbf{Q}_{n(m+1)}$, and translate (\ref{minT_PB}) with respect to such matrices:
\begin{itemize}
  \item {\em Objective function}.  We have
\eqn && \hspace{-0.8cm}\tr\int_{-\pi}^{\pi} \Psi_y(e^{i\vartheta})\frac{\de \vartheta}{2\pi}=\int \tr(\Delta(e^{i\vartheta})Y \Delta(e^{i\vartheta})^*)\frac{\de \vartheta}{2\pi}\nn\\ &&\hspace{0.0cm}=\tr\left (Y \int_{-\pi}^\pi \Delta(e^{i\vartheta})^*\Delta(e^{i\vartheta})\frac{\de\vartheta}{2\pi} \right)=\tr(Y)\nn\eeqn
where we exploited the fact that
\eq \int_{-\pi}^\pi e^{ik\vartheta}\frac{\de \vartheta}{2\pi}=\left\{
                                                   \begin{array}{ll}
                                                     1, & k=0 \\
                                                     0, & k\neq 0.
                                                   \end{array}
                                                 \right.
\eeq

\item {\em Conditions $\Psi_y,\Psi_z\succeq 0$}. The condition $Y\succeq 0$ implies that $\Psi_y(e^{i\vartheta})=\Delta (e^{i\vartheta})Y \Delta(e^{i\vartheta})^*\succeq 0$ for each $\vartheta\in[-\pi,\pi]$. On the other hand if $\Psi_y\succeq 0$,
there exists $\Gamma(e^{i\vartheta})=\sum_{k=0}^m e^{-ik \vartheta} C_k$, with $C_k\in\Rs^{n \times l}$, such that $\Psi_y(e^{i\vartheta})=\Gamma(e^{i\vartheta})\Gamma(e^{i\vartheta})^*$.
Hence, $\Psi_y(e^{i\vartheta})=\Delta(e^{i\vartheta}) Y\Delta(e^{i\vartheta})^*$ with
\eq Y=\left[
        \begin{array}{c}
          C_0 \\
          C_1 \\
          \vdots \\
          C_m \\
        \end{array}
      \right]\left[
               \begin{array}{cccc}
                 C_0^T & C_1^T & \ldots & C_m^T \\
               \end{array}
             \right]
\eeq which is positive semidefinite. Thus, we can replace $\Psi_y\succeq 0$ with $Y\succeq 0$, and similarly $\Psi_z\succeq 0$ with $Z\succeq 0$.
\item   {\em Condition $\Psi_x=\Psi_y+\Psi_z$}. Let $\Psi_x(e^{i\vartheta})=\sum_{k=-m}^m e^{-ik \vartheta}R_k$, $\Psi_y(e^{i\vartheta})=\sum_{k=-m}^m e^{-ik \vartheta}P_k$, $\Psi_z(e^{i\vartheta})=\sum_{k=-m}^m e^{-ik \vartheta}Q_k$ with $R_k=R_{-k}^T$, $P_k=P_{-k}^T$, $Q_k=Q_{-k}^T$. Thus, the equality constraint may be rewritten as $\sum_{k=-m}^m e^{-ik \vartheta}(P_k+Q_k-R_k)=0$ which is equivalent to
    \eq \label{cond_coi_coefficienti} P_k+Q_k=R_k, \; k=0\ldots m.\eeq
    Note that $R_k$s, $P_k$s and $Q_k$s are the coefficients of the pseudo-polynomial matrices $\Psi_x$, $\Psi_y$ and $\Psi_z$, respectively.
    In view of (\ref{representation_matricial}), the $k$-th coefficients of $ \Psi_y$ and $\Psi_z$ are given by $\mathbf{D}_k(Y)$
and $\mathbf{D}_k(Z)$, respectively. Thus, (\ref{cond_coi_coefficienti}) is equivalent to
\eq \mathbf{D}_k(Y)+\mathbf{D}_k(Z)=R_k,\; k=0\ldots m.\eeq Finally, by exploiting the linearity of $\mathbf{D}_k$s, we obtain
\eq \mathbf{D}_k(Y+Z)=R_k,\; k=0\ldots m.\eeq
\item {\em Condition $\Psi_z$ diagonal}. By exploiting argumentations similar to the ones of the previous point, we get that the condition is equivalent to
\eq \mathbf{D}_k(Z) \hbox{ diagonal}, \; k=0\ldots m.\eeq
\end{itemize}
Hence, Problem (\ref{minT_PB}) is equivalent to
\eqn \label{PB_matricial} &\underset{Y,Z\in\mathbf{Q}_{n(m+1)}}{\min}& \tr(Y)\nn\\
   &\hbox{ subject to }&   Y,Z\succeq 0\nn\\
   && \mathrm{D}_k(Z) \hbox{ diagonal } k=0\ldots m\nn\\
  &&   \mathrm{D}_k(Y+Z)=R_k\; k=0\ldots m \eeqn
and a solution to (\ref{minT_PB}) is given by $\hat \Psi_y(e^{i\vartheta})=\Delta(e^{i\vartheta}) \hat Y \Delta(e^{i\vartheta})^*$ and $\hat \Psi_z(e^{i\vartheta})=\Delta(e^{i\vartheta}) \hat Z \Delta(e^{i\vartheta})^*$, where $(\hat Y,\hat Z)$ is solution to (\ref{PB_matricial}).
We conclude that a solution to (\ref{minT_PB}) may be easily computed by solving (\ref{PB_matricial}).
\smallskip \smallskip \rem In the case that $x$, $y$, $z$ are MA processes of order $m_x$, $m_y$ and $m_z$, respectively, and such that $m_x<m_y\leq m_z$, we define $\Psi_x, \Psi_y,\Psi_z\in \Qc_{n,m_z}$ and (\ref{PB_matricial}) becomes
\eqn  &\underset{Y,Z\in\mathbf{Q}_{n(m_z+1)}}{\min}& \tr(Y)\nn\\
   &\hbox{ subject to }&   Y,Z\succeq 0\nn\\
   && \mathrm{D}_k(Z) \hbox{ diagonal } k=0\ldots m_z\nn\\
  &&   \mathrm{D}_k(Y+Z)=R_k\; k=0\ldots m_x \nn \\
  && \mathrm{D}_k(Y+Z)=0\; k=m_x+1\ldots m_y\nn\\
  && \mathrm{D}_k(Y)=0\; k=m_y+1\ldots m_z \eeqn
  and it is not difficult to show that the previous results, in particular Proposition \ref{prop_esistenza}, still holds.\erem

\section{Simulation studies}\label{section_sim}

\subsection{Performance of the relaxed problem} \label{subsection_performance_relaxed}
We start by testing the tightness of the convex relaxation (\ref{minT_PB}). We consider 10 dynamic factor models (\ref{DFA_model}) with $n=10$ and $m=5$ whose coefficients $A_k$s $B_k$s are randomly generated. These models differ by the number of common factors, i.e. $r$. For each model we solve problem (\ref{minT_PB}) by using the true spectral density $\Psi_x$ of the observed process.
Let $\Psi_y$ be the true spectral density of $y$, and $\hat \Psi_y$ the estimate provided by  (\ref{minT_PB}). We compute the relative estimation error, averaged on the unit circle, of $\hat \Psi_y$:
\eq \mathrm{m}_{\mathrm{e}_{\Psi_y}}=\int_{-\pi}^{\pi}\frac{\| \Psi_y(e^{i\vartheta})-\hat \Psi_y(e^{i\vartheta})\|_2}{\| \Psi_y(e^{i\vartheta})\|_2}\frac{\de\vartheta}{2\pi}.\eeq
  In the following table the relative error $\mathrm{m}_{\mathrm{e}_{\Psi_y}}$ for the 10 different models is shown.
\smallskip \smallskip \begin{center}
\begin{tabular}{|c|c|c|c|}
  \hline
  $r$ & $\mathrm{m}_{\mathrm{e}_{\Psi_y}}$ & $r$ & $\mathrm{m}_{\mathrm{e}_{\Psi_y}}$ \\
  \hline
  1 & $8.72\times 10^{-10}$  & 6 & $ 5.02\times 10^{-5}$\\
  2 & $4.56\times 10^{-10}$ &  7 & $ 2.91\times 10^{-2}$\\
  3 & $1.76\times 10^{-9}$ &   8 & $ 3.15\times 10^{-2}$\\
  4 & $ 1.25 \times 10^{-9}$ & 9 & $ 5.51 \times 10^{-2}$\\
  5 & $ 1.63 \times 10^{-9}$ & 10 & $8.25 \times 10^{-2}$\\
 \hline
\end{tabular}
\end{center} \smallskip \smallskip
The estimated decomposition is exact when the number of common factors is small, $r\leq 5$, whereas
it is not for a large number of common factors, $r\geq 6$. At this point it is worth recalling that the decomposition (\ref{DFA_freq_repr})
is generically unique, and thus identifiable from $\Psi_x$, for $r\leq n- \sqrt{n}$, \cite{deistler2007}. In our case $n- \sqrt{n} \cong 6.84$. We conclude that the relaxed formulation is able to recover all the identifiable decompositions except for the case $r=6$.

\subsection{Factor Analysis}
\begin{figure}[thpb]
      \centering
      \includegraphics[width=9.2cm]{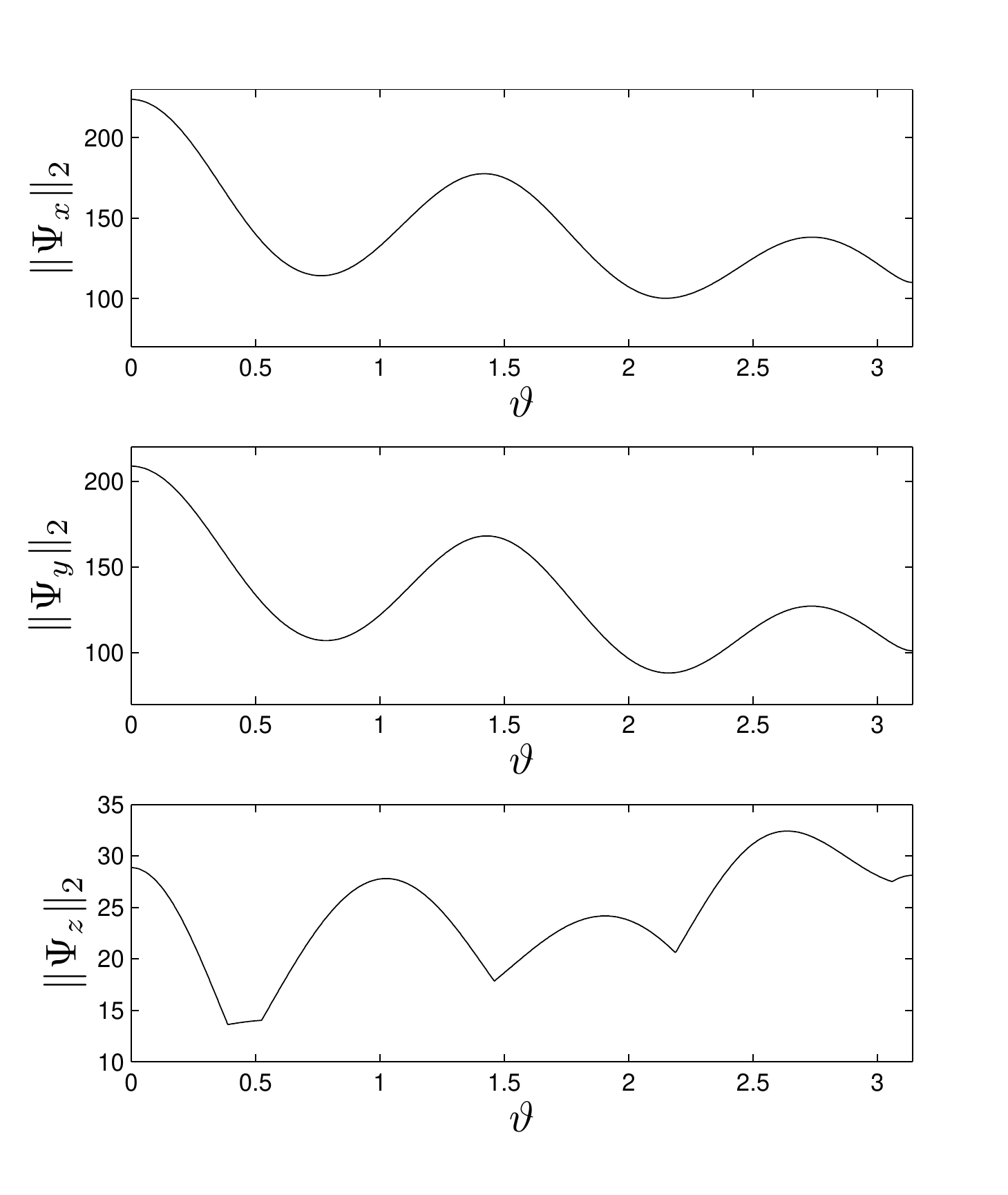}
      \caption{Spectral norm of $\Psi_x$, $\Psi_y$ and $\Psi_z$ of the dynamic factor model with $r=3$.}
      \label{fig_factor_model}
   \end{figure}

\begin{figure}[thpb]
      \centering
      \includegraphics[width=9.2cm]{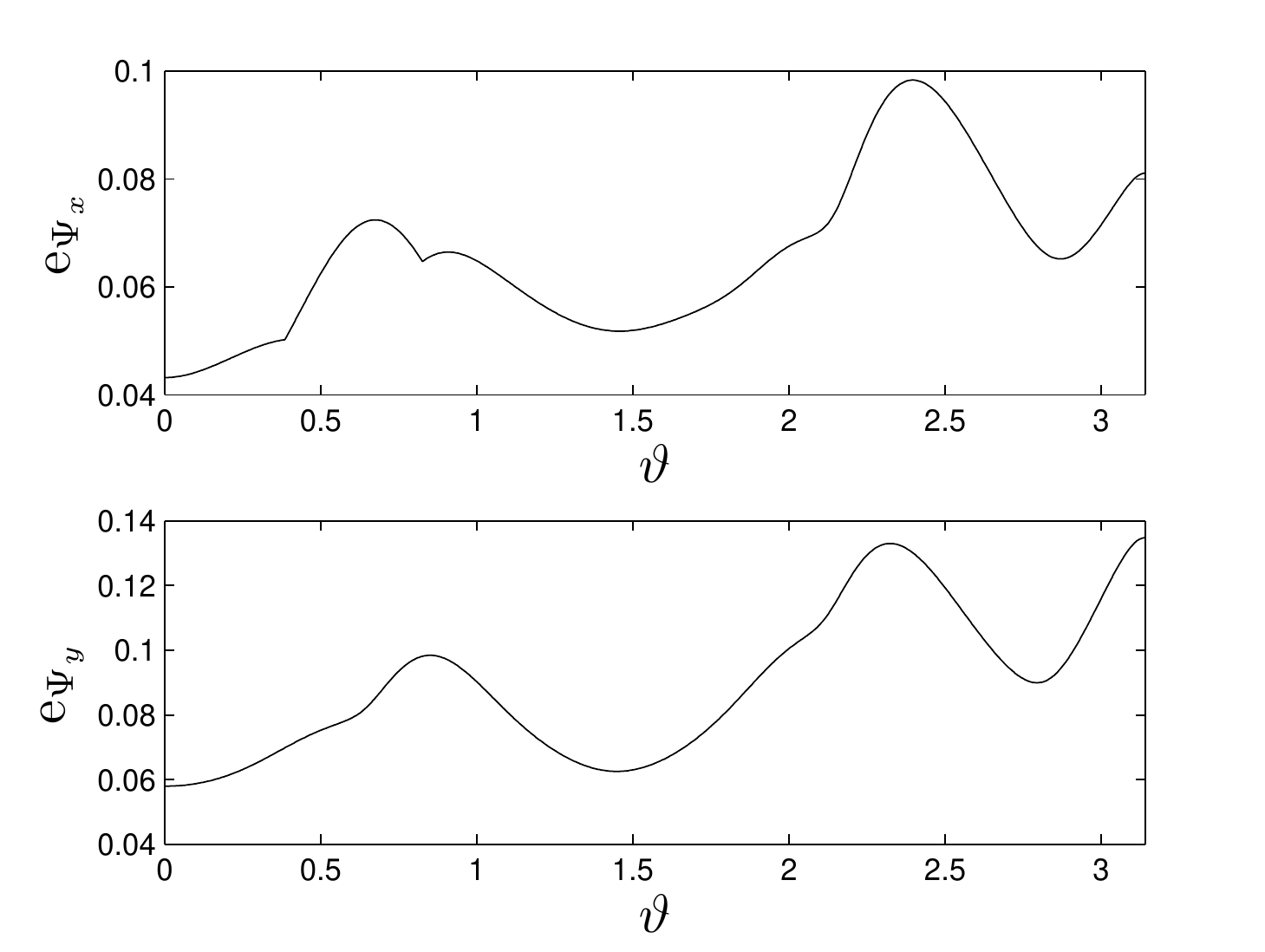}
      \caption{Relative estimation error of $\Psi_x$ and $\Psi_y$.}
      \label{fig_relative_error}
   \end{figure}

\begin{figure}[thpb]
      \centering
      \includegraphics[width=9cm]{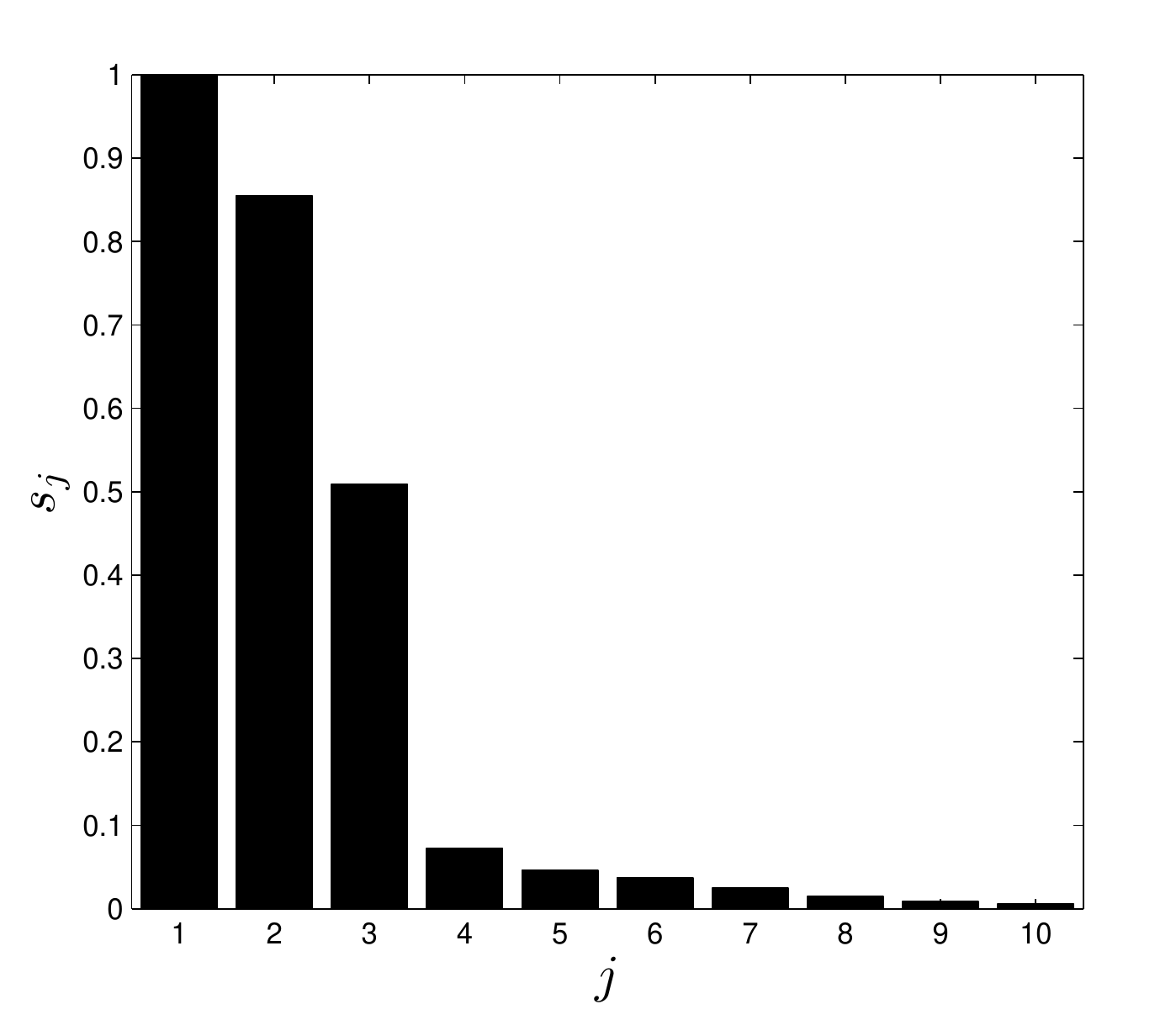}
      \caption{Normalized singular values over the unit circle of $\hat \Psi_y$.}
      \label{fig_eigenvalues}
   \end{figure}

We consider the dynamic factor model of Section \ref{subsection_performance_relaxed} with $r=3$.
In Figure \ref{fig_factor_model}
  we depict the spectral norm of $\Psi_x$, $\Psi_y$ and $\Psi_z$ at each frequency. From (\ref{DFA_model}) we generate $N=6000$ samples of $x$: $\mathrm{x}(1),\mathrm{x}(2),\ldots \mathrm{x}(6000)$.
   We apply then the identification procedure of Section \ref{section_ID_procedure} for characterizing common and specific factors from the data.
   We denote by $\hat \Psi_x$ the estimate of $\Psi_x$ computed by Durbin's method, and $\hat \Psi_y$ the estimate of $\Psi_y$ obtained by solving Problem (\ref{minT_PB}).
   We define the relative estimation errors at each frequency:
   \eqn  \mathrm{e}_{\Psi_x}(e^{i\vartheta})&=&\frac{\| \Psi_x(e^{i\vartheta})-\hat \Psi_x(e^{i\vartheta})\|_2}{\|\Psi_x(e^{i\vartheta})\|_2}\nn\\
   \mathrm{e}_{\Psi_y}(e^{i\vartheta})&=&\frac{\| \Psi_y(e^{i\vartheta})-\hat \Psi_y(e^{i\vartheta})\|_2}{\|\Psi_y(e^{i\vartheta})\|_2}\eeqn
The errors graph is displayed in Figure \ref{fig_relative_error}.  We note that $e_{\Psi_x}$ and $e_{\Psi_y}$ take similar
   values for $\vartheta\in[-\pi,\pi]$. This means that the estimation error is mainly imputable
   to Durbin's method for estimating $\Psi_x$ from the data.
Finally, we define \eq s_j:= \underset{\vartheta\in[-\pi,\pi]}{\max} \frac{\sigma_j(\hat \Psi(e^{i\vartheta}))}{\sigma_1(\hat \Psi(e^{i\vartheta}))}\eeq
where $\sigma_j(\hat \Psi(e^{i\vartheta}))$ is the $j$-th largest singular value of $\hat \Psi_y$ at $\vartheta$.
Hence, $s_j$ can be understood as the $j$-th largest normalized singular value over the unit circle of $\hat \Psi_y$. Those quantities are plotted in Figure \ref{fig_eigenvalues}.   The plot suggests  that we can safely approximate $\mathrm{rank}(\hat \Psi_y)=3$. Accordingly, we recover the exact number of common factors. Finally, we obtained similar results with different samples and by
changing the factor model.

\section{Conclusion}\label{section_conclusion}
In this paper we proposed an identification procedure for factor analysis of MA processes. Here, the challenging step is to solve a minimum rank problem.
We proposed a convex optimization problem approximating the NP-hard problem. Simulation studies point out that the convex problem is able to recover a correct solution
in most of the cases. Finally, we tested the identification procedure: simulations show this method is able to identify common and specific factors with satisfactory accuracy.

\end{document}